\documentclass[12pt,a4paper]{article}
\usepackage{mathrsfs}

\usepackage{amssymb}

\usepackage{amsmath}

\newtheorem{theorem}{Theorem}[section]
\newtheorem{definition}[theorem]{Definition}
\newtheorem{lemma}[theorem]{Lemma}

\begin{document}

\title{GR\"{O}BNER-SHIRSHOV BASES: SOME NEW RESULTS\footnote{Supported by the
NNSF of China (No.10771077) and the NSF of Guangdong Province
(No.06025062).}}

\author{L. A. Bokut\footnote {Supported by the RFBR  and
the Integration Grant of the SB RAS (No. 1.9).}\\
{\small School of Mathematical Sciences, South China Normal
University,}\\
{\small Guangzhou 510631,  China}\\
{\small Sobolev Institute of Mathematics, Russian Academy of
Sciences}\\
{\small Siberian Branch, Novosibirsk 630090, Russia}\\
{\small bokut@math.nsc.ru}\\
\\
 Yuqun Chen\footnote {Corresponding author.}\\
{\small School of Mathematical Sciences, South China Normal
University,}
\\ {\small Guangzhou 510631,  China}\\
{\small  yqchen@scnu.edu.cn}}

\date{}
 \maketitle

\maketitle \noindent\textbf{Abstract:} In this survey article, we
report some new results of Gr\"{o}bner-Shirshov bases, including new
Composition-Diamond lemmas,  applications of some known
Composition-Diamond lemmas and content of some expository papers.

\maketitle \noindent\textbf{Key words: } Composition-Diamond lemma;
Group; HNN-extension; Schreier extension; Dialgebra; Lie algebra;
Module; Chinese monoid.

 \section{Introduction}
In this survey, we report the activities of the first author who has
been staying in the South China Normal University, at Guangzhou, in
Spring 2006, Spring 2007 and November-December 2007. With the
participation of the second author, we are running an algebra
seminar with 12 students (most of them are master degree students)
and some young teachers at 5 times per week, two hours each session.
The subjects of this seminar include from Combinatorial group
theory, Free Lie algebras, Semi-simple Lie algebras to
Non-associative algebras, Conformal algebras, Quantum groups,
Semigroups and Dialgebras with emphasizing in Gr\"{o}bner-Shirshov
bases.

The second author visited Sobolev Institute of Mathematics at
Novosibirsk as a visiting professor in July-October, 2006.

As the result of all these activities, more than 10 papers have been
prepared. We now give a brief survey for some of the papers. We also
mention some papers which were done during the first author's
visiting to Prof. K. P. Shum, the Chinese University of Hong Kong,
2003-2005.

We divide these papers into three blocks:
\begin{enumerate}
\item[(i)] \textit{New Composition-Diamond (CD-) lemmas,}
\item[(ii)] \textit{Applications of known CD-lemmas,}
\item[(iii)] \textit{Expository papers.}
\end{enumerate}

We first explain  what it means of ``CD-Lemma"  for a class (variety
or category) $\cal{M}$ of linear $\Omega$-algebras over a field $k$
(here $\Omega$ is a set of linear operations on $\cal{M}$) with free
objects.

\textbf{$\cal{M}$-CD-Lemma} Let $\cal{M}$ be a class of (in general,
non-associative) $\Omega$-algebras, $Free_{\cal{M}}(X)$ a free
$\Omega$-algebra in $\cal{M}$ generated by $X$ with a linear base
consisting of ``normal (non-associative $\Omega$-) words" $[u]$, $S
\subset \ Free_{\cal{M}}(X)$ a subset and $<$ a monomial well order
on normal
 words. Let $S$ be a Gr\"{o}bner-Shirshov basis (this
means that any ``composition" of elements in $S$ is ``trivial").
Then

\begin{enumerate}
\item[(a)] If $f\in Ideal(S)$, then $[\bar{f}]=[a\bar{s}b]$, where $[\bar{f}]$
is the ``leading monomial" of $f$ and $[asb]$ is a ``normal
$s$-word", $s\in S$.

\item[(b)] $Irr(S)$ =$\{[u]| [u] \neq[a\bar{s}b], s\in S,
[asb] \ is \ a \ normal \ s-word\}$ is a linear basis of the algebra
${\cal M}(X|S)$ with defining relations $S$.
\end{enumerate}

In many cases, each of conditions (a) and (b) is equivalent to the
condition that $S$ is a Gr\"{o}bner-Shirshov basis in
$Free_{\cal{M}}(X)$.  But in some of our ``new CD-Lemmas", this is
not the case.
Here, $\cal{M}$ may be as follows:\\

--Associative algebras (\cite{Sh}, \cite{b76}, \cite{b}).

A free associative algebra is $k\langle X\rangle$, the algebra of
non-commutative polynomials; normal words are words on $X$; the A.
I. Shirshov's composition $(f,g)_w$ is equal to $fb-ag$, if
$w=\bar{f}b=a\bar{g}$, deg$(\bar{f})+$deg$(\bar{g})>$deg$(w)$, or
$f-agb$, if $w=\bar{f}=a\bar{g}b$, where $a,b\in X^*$ and $X^*$ the
free monoid generated by $X$; a polynomial $h$ is called trivial
$mod(S)$ if it goes to 0 by using the Eliminations of Leading Words
(ELW) of $S$ (see below an equivalent definition).

\ \ --Lie algebras (\cite{Sh}).

A free Lie algebra is $Lie(X)$, the algebra of Lie polynomials in
$k\langle X\rangle$ (this theorem was proved by W. Magnus and E.
Witt); by the normal words we mean the non-associative
Lyndon-Shirshov words $[u]$ on $X$; the leading word $\bar{f}$ of a
Lie polynomial $f$ is the same as the associative polynomial; A. I.
Shirshov's composition $[f,g]_w$ of two Lie polynomials is its
associative composition with some extra bracketing defined in
\cite{S58}; a normal $s$-word for $s\in Lie(X)$ has the form $[asb]$
with extra bracketing as before.

\ \ --Commutative algebras (\cite{bu65}, \cite{bu70}).

A free commutative algebra is $k[X]$, the algebra of polynomials  on
$X$ over a field $k$;  normal words are monomials; the composition
$S(-,-)$ is the operation of taking the  B. Buchberger's
$S$-polynomial: $S(f,g)=fb-ag$ for any polynomials $f,g$, where
$w=\bar{f}b=a\bar{g}=l.c.m(\bar{f},\bar{g})$ and
deg$(\bar{f})+$deg$(\bar{g})>$deg$(w)$.

\ \

--(Commutative, anti-commutative) non-associative algebras
(\cite{S3}).

There is only composition of inclusion in the cases.

\ \

--Lie superalgebras (\cite{Mik1}).

The Composition-Diamond lemma for Lie superalgebras is known and
proved.

\ \

--Grassmann algebras (\cite{stokes}).

There is new composition of multiplication by a monomial.

\ \

--Supercommutative associative superalgebras (\cite{Mik1},
\cite{Mik2}).

There is new composition of multiplication by a monomial.

\ \

--Conformal associative algebras $(C, \ (n), \ n\geq 0, \ D)$
(\cite{bfk}).

There are 6 types of compositions including inclusion, intersection,
$D$-inclusion, $D$-intersection, left (right) multiplication by a
generator. The condition (a) ((b)) in the CD-lemma is not equivalent
to the condition that $S$ is a Gr\"{o}bner-Shirshov basis.

\ \

--Modules (\cite{kl1}, \cite{ch04}).

\section{CD-lemma for associative algebras}

In this section, we cite some concepts and results from the
literature which are related to the Gr\"{o}bner-Shirshov bases for
the associative algebras.

\begin{definition} (\cite{Sh}, see also \cite{b72}, \cite{b76}) \
Let $f$ and $g$ be two monic polynomials in \textmd{k}$\langle
X\rangle$ and $<$ a well order on $X^*$. Then, there are two kinds
of compositions:

$(i)$ If \ $w$ is a word such that $w=\bar{f}b=a\bar{g}$ for some
$a,b\in X^*$ with deg$(\bar{f})+$deg$(\bar{g})>$deg$(w)$, then the
polynomial
 $(f,g)_w=fb-ag$ is called the intersection composition of $f$ and
$g$ with respect to $w$.

$(ii)$ If  $w=\bar{f}=a\bar{g}b$ for some $a,b\in X^*$, then the
polynomial $(f,g)_w=f - agb$ is called the inclusion composition of
$f$ and $g$ with respect to $w$.

\end{definition}

\begin{definition}(\cite{b72}, \cite{b76}, \cite{Sh})
Let $S\subset k\langle X\rangle$ with each $s\in S$ monic. Then the
composition $(f,g)_w$ is called trivial modulo $(S,w)$ if
$(f,g)_w=\sum\alpha_i a_i s_i b_i$, where each $\alpha_i\in k$,
$a_i,b_i\in X^{*}$, $s_i\in S$ and $\overline{a_i s_i b_i}<w$. If
this is the case, then we write $ (f,g)_w\equiv0\quad mod(S,w) $
\end{definition}

\begin{definition} (\cite{b72}], \cite{b76},  \cite{Sh}) \
We call the set $S$ with respect to the well order $``<"$ a
Gr\"{o}bner-Shirshov set (basis) in $k\langle X\rangle$ if any
composition of polynomials in $S$ is trivial modulo $S$.
\end{definition}

If a subset $S$ of $k\langle X\rangle$ is not a Gr\"{o}bner-Shirshov
basis, then we can add to $S$ all nontrivial compositions of
polynomials of $S$, and by continuing this process (maybe
infinitely) many times, we eventually obtain a Gr\"{o}bner-Shirshov
basis $S^{comp}$. Such a process is called the Shirshov algorithm.
It is an infinite algorithm as well as Kruth-Bendix algorithm (see
\cite{knuth}).

A well order $>$ on $X^*$ is monomial if it is compatible with the
multiplication of words, that is, for $u, v\in X^*$, we have
$$
u > v \Rightarrow w_{1}uw_{2} > w_{1}vw_{2},  \ for \  all \
 w_{1}, \ w_{2}\in  X^*.
$$
A standard example of monomial order on $X^*$ is the deg-lex order
to compare two words first by degree and then lexicographically,
where $X$ is a linearly ordered set.

The following lemma was proved by Shirshov \cite{Sh} for the free
Lie algebras (with deg-lex ordering) in 1962 (see also Bokut
\cite{b72}). In 1976, Bokut \cite{b76} specialized the approach of
Shirshov to associative algebras, see also Bergman \cite{b}. For
commutative polynomials, this lemma is known as the Buchberger's
Theorem in \cite{bu65} and \cite{bu70}.

\begin{lemma}\label{l1}
(Composition-Diamond Lemma) \ Let $k$ be a field, $A=k \langle
X|S\rangle=k\langle X\rangle/Id(S)$ and $<$ a monomial order on
$X^*$, where $Id(S)$ is the ideal of $k \langle X\rangle$ generated
by $S$. Then the following statements are equivalent:
\begin{enumerate}
\item[(i)] $S $ is a Gr\"{o}bner-Shirshov basis in $k \langle
X\rangle$.
\item[(ii)] $f\in Id(S)\Rightarrow  \bar{f}=a\bar{s}b$
for some $s\in S$ and $a,b\in  X^*$.
\item[(iii)] $Irr(S) = \{ u \in X^* |  u \neq a\bar{s}b ,s\in S,a ,b \in X^*\}$
is a basis of the algebra $A=k\langle X | S \rangle$.
\end{enumerate}
\end{lemma}

\section{New CD-lemmas}
\subsection{CD-Lemma and HNN-extensions}

Y. Q. Chen and C. Y. Zhong in \cite{cz1} give a version of CD-lemma
in which the order may not be monomial. A Gr\"{o}bner-Shirshov basis
for HNN extensions of groups is obtained by using the new CD-lemma.
This is the first paper to give a Gr\"{o}bner-Shirshov basis by
using a non-monomial order.

\begin{theorem}\label{3.1}
(\cite{cz1}) \ Let $S\subseteq k\langle X\rangle$ and $``<"$ a well
order on $X^*$ such that
\begin{enumerate}
\item[(i)] \ $\overline{asb}=a\bar{s}b$ for any $a,b\in X^*, \ s\in S$;
\item[(ii)] \ for each composition $(s_1,s_2)_w$ in $S$, there exists a presentation
$$
(s_1,s_2)_w=\sum_{i}\alpha_{i}a_it_ib_i, \ a_i\bar{t_i}b_i<w, \ \
\mbox{ where } \ t_i\in S, \ a_i,b_i\in X^*, \ \alpha_{i}\in k
$$
such that for any $c,d\in X^*$, we have $ca_i\bar{t_i}b_id<cwd$.
\end{enumerate}
Then, the following statements hold.
\begin{enumerate}
\item[(i)] \ $S$ is a Gr\"{o}bner-Shirshov basis in $k \langle
X\rangle$.
\item[(ii)] \ For any $f\in k\langle X\rangle, \
f\in Id(S)\Rightarrow \bar{f}=a\bar{s}b$ for some $s\in S, \ a,b\in
X^*$.
\item[(iii)] \
The set
$$
Irr(S)=\{u\in X^*|u\neq a\bar{s}b,s\in S,a,b\in X^*\}
$$
is a linear basis of the algebra $k\langle X|S\rangle$.
\end{enumerate}
\end{theorem}

We call the order satisfying the conditions in Theorem \ref{3.1} an
$S$-weak monomial order.

Let $G=gp\langle H,t|t^{-1}at=\varphi(a), a\in A \rangle$ be an
HNN-extension of a group H, where $A$ is a subgroup of $H$ and
$\varphi$ a group isomorphism.  By using Theorem \ref{3.1}, it is
proved in [\cite{cz1}] that there exists an explicit
Gr\"{o}bner-Shirshov basis $S$ of $G$ relative to some explicit
$S$-weak monomial order such that the set $Irr(S)$ of
$S$-irreducible words coincides with the set of normal forms in the
Normal Form Theorem for HNN-extensions (see [\cite{b66-7}] and
[\cite{ly}]).

\subsection{ Dialgebreas }

In this section, we report some recent results of L. A. Bokut, Y. Q.
Chen and C. H. Liu  \cite{bcl}.

Let $D(X)$ be a free dialgebra (J.-L. Loday, 1995, \cite{Lo95}),
where multiplications $``\vdash",\ \ ``\dashv"$ are both associative
and for any $a,b,c\in D(X)$,
\begin{eqnarray*}
a\dashv(b\vdash c)&=&a\dashv b\dashv c, \\
(a\dashv b)\vdash c&=&a\vdash b\vdash c, \\
a\vdash(b\dashv c)&=&(a\vdash b)\dashv c.
\end{eqnarray*}

A linear basis of $D(X)$ consists of normal diwords
$$
[u]= x_{-m}\vdash \cdots \vdash x_{0}\dashv\cdots \dashv x_k=
x_{-m}\cdots \dot{x_{0}}\cdots x_k,
$$
where $x_i\in X, m, k\geq0, x_0$ is the center of $[u]$ (see J.-L.
Loday, 1995, \cite{Lo95}). We define the deg-lex order $[u]<[v]$, by
using the lex-order of the weight $wt [u]=(k+m+1, m,x_{-m},\dots,
x_k)$.

Now, for $f, g\in S$, we define the compositions of inclusion,
intersection and left(right) multiplication by a letter.

We call the set $S$ a Gr\"{o}bner-Shirshov set (basis) in $D(X)$ if
any composition of polynomials in $S$ is trivial modulo $S$ (and
$[w])$.

\begin{theorem}  (\cite{bcl} CD-Lemma for dialgebras) Let $S\subset D(X)$ be a monic set and
the order $<$ as before. Then $(i)\Rightarrow (ii)\Leftrightarrow
(iii)\Rightarrow (iv)$, where
\begin{enumerate}
\item[(i)] \ $S$ is a Gr\"{o}bner-Shirshov basis in $D(X)$.
\item[(ii)] \
For any $f\in D(X), \  f\in Id(S)\Rightarrow
[\overline{f}]=[a[\overline{s}]b]$ for some $s\in S, \ a,b\in [X^*]$
 and $[asb]$ a normal $S$-diword.
\item[(iii)] \
The set $
Irr(S)=\{u\in [X^*]|u\neq [a[\overline{s}]b],s\in S,a,b\in [X^*],\\
\quad \quad\quad \quad\quad \quad\quad \quad [asb] \mbox{ is normal
S-diword}\} $ is a linear basis of the dialgebra $D(X|S)$.
\item[(iv)] \
 Each composition is trivial modulo $S$.
\end{enumerate}
\end{theorem}

As an application of the above theorem, we obtain a
Gr\"{o}bner-Shirshov basis for the universal enveloping algebra of a
Leibniz algebra. It is the PBW-Theorem for the  Leibniz algebras.
This is the third proof of the theorem after M. Aymon and P. P.
Grival (2003) in \cite{ag}, and P. Kolesnikov (2007) in \cite{ko}.

Recall that a Leibniz algebra $L$ is a non-associative algebra with
a multiplication $[xy]\in L$ such that $[[xy]z]-[[xz]y]-[[yz]x]=0$
(see \cite{Lo95}). For any dialgebra $(D,\dashv,\vdash)$, the linear
space $D$ with the multiplication $[xy]=x\dashv y-y\vdash x$ is  a
Leibniz algebra. For any Leibniz algebra
$L=Lei(\{e_i\}_I|[e_ie_j]=\sum_k\alpha_{ij}^ke_k, \ i,j\in I)$, one
can define the universal enveloping $D$-algebra
$U(L)=D(\{e_i\}_I|e_i\dashv e_j-e_j\vdash
e_i=\sum_k\alpha_{ij}^ke_k, \ i,j\in I)$, where $\{e_i\}_I$ is a
basis of $L$.

\begin{theorem}(\cite{bcl})
Let $\mathcal{L}$ be a Leibniz algebra over a field $k$ with the
product $\{,\}$. Let $\mathcal{L}_0$ be the subspace of
$\mathcal{L}$ generated by the set $\{\{a,a\}, \{a,b\}+\{b,a\} \ | \
a,b\in \mathcal{L}\}$. Let $\{x_i|i\in I_0\}$ be a basis of
$\mathcal{L}_0$ and $X=\{x_i|i\in I\}$ a linearly ordered basis of
$\mathcal{L}$ such that $I_0\subseteq I$. Let $(D(X),\dashv,\vdash)
$ be the free dialgebra and the order $<$ on $[X^*]$ as before. Let
$S$ be the set which consists of the following polynomials:
\begin{eqnarray*}
1.&&f_{ji}=x_j\vdash x_i-x_i\dashv x_j+\{x_i,x_j\} \ \  \ \ \  \ \ \
\ \ \  \ \ \  \ \ \  \ \ \ \
 \  \  \ \ \  (i,j\in I)\\
2.&&f_{ji\vdash t}=x_j\vdash x_i\vdash x_t-x_i\vdash x_j \vdash
x_t+\{x_i,x_j\}\vdash x_t \ \ \ \ \  \    (i,j,t\in I, \ j>i)\\
3.&&h_{i_0\vdash t}=x_{i_0}\vdash x_t \ \ \ \ \ \ \ \  \  \ \ \ \  \
\  \ \ \ \ \ \ \ \ \ \ \  \  \ \ \ \  \  \  \ \ \ \ \ \ \ \ \ \ \  \
\ \
(i_0\in I_0, \ t\in I)\\
4.&&f_{t\dashv ji}=x_t\dashv x_j\dashv x_i-x_t\dashv x_i\dashv
x_j+x_t\dashv \{x_i,x_j \} \ \ \ \
 \ \  (i,j,t\in I, \ j>i)\\
5.&&h_{t\dashv i_0}=x_t\dashv x_{i_0} \ \ \ \ \ \ \ \  \ \ \ \ \  \
\  \ \ \ \ \ \ \ \ \ \ \  \  \ \ \ \  \  \  \ \ \ \ \ \ \ \ \ \ \  \
\ \   (i_0\in I_0, \ t\in I).
\end{eqnarray*}
Then
\begin{enumerate}
\item[(i)] \ $S$ is a Gr\"{o}bner-Shirshov basis in $(D(X)$.
\item[(ii)] \ The set
$$
\{x_j\dashv x_{i_1}\dashv \dots \dashv x_{i_k} \ | \ j\in I, i_p\in
I-I_0, \ 1\leq p\leq k, \ i_1\leq \dots \leq i_k, \ k\geq 0\}
$$
is a linear basis of the universal enveloping algebra
$U(\mathcal{L})=D(X|S)$. In particular, $\mathcal{L}$ can be
embedded into $U(\mathcal{L})$.
\end{enumerate}
\end{theorem}

\subsection{ Free $\Gamma$-algebras $k\langle
X;\Gamma\rangle$}

In this section, we summary the results given by  L. A. Bokut and K.
P. Shum  \cite{bs1}.

Let $X$  be a set, $\Gamma$ a group, $\Gamma(x)$, $\Gamma'(x)$
isomorphic subgroups, $x\in X$. Then the algebra $k\langle
X;\Gamma\rangle$ with defining relations
$$
\gamma x=x\gamma' \ (\gamma \in \Gamma (x),\ \gamma'\in \Gamma'(x),
x\in X), \ \ \ \gamma\delta =\mu \ (\gamma,\delta,\mu\in \Gamma)
$$
is called $free$ $\Gamma$-$algebra$.

 A linear basis of $k\langle X;\Gamma\rangle$
consists of $\Gamma$-$words$
$$
u=\gamma_0x_{i_1}\gamma_{1}\cdots x_{i_k}\gamma_k,\ x_i\in X,
\gamma_i\in \Gamma, k\geq0,
$$
which are equivalent under transformations $\gamma x\rightarrow
x\gamma'$ above.

 We input a quasi-order on $\Gamma$-words:
$$
u\leq v \Leftrightarrow  [u]\leq [v],
$$
where $[u]=x_{i1}\cdots x_{ik}$ is the projection of $u$, and $ [u]
\leq [v]$ a monomial order on $X^*$.

 A $\Gamma$-polynomial $f$ may have  several leading
monomials of
 $\bar{f}$. We call $f$ a $strong$ polynomial if  $\bar{f}$ is
 unique.
  We define compositions of inclusion and
intersection of two strong $\Gamma$-polynomials, and a strong
$\Gamma$-Gr\"{o}bner-Shirshov basis. The later is a set of strong
$\Gamma$-polynomials that is closed under compositions.

\begin{theorem}   Let  $k\langle X;\Gamma\rangle$ be a  free strong
$\Gamma$-algebra, $S\subset
  k\langle X;\Gamma\rangle$ a  strong $\Gamma$-Gr\"{o}bner-Shirshov basis. Then
\begin{enumerate}
  \item[(a)] If $f\in Id(S)$, then $\bar{f}=a\overline{s}b$, where
$\bar{f}$ is a leading monomial of $f$, $s\in S$, $a,b\
\Gamma$-words.

 \item[(b)] $Irr(S)=\{u\neq a\bar{s}b|s\in S, a, b \ $ are \ $
  \Gamma$-words\} is a linear basis of $k\langle X;\Gamma|S\rangle$.
  \end{enumerate}
\end{theorem}

There are many examples of $\Gamma$-algebras with strong
$\Gamma$-Gr\"{o}bner-Shirshov bases.

  \textbf{(a)} \ \ Group algebras of universal groups  $G(R^*)$ of multiplicative
  semigroups $R^*$ of some rings $R$.

Let $R= \overline{k(S)}$, where $S=sgp\langle X | w_ih_i=u_if_i,
w_i,h_i,u_i,f_i\in X\rangle$, $k$ a field, $k(S)$ the semigroup
algebra, $ \overline{k(S)}$ the algebra of formal series over $S$.
In particular, if $S$ is a free semigroup, then $ \overline{k(S)}=
\overline{k\langle X\rangle}$ is the Magnus algebra of formal series
over $X$.

 These examples are from Bokut's solution to the Malcev
embedding problem: There exists a semigroup $S$ such that
$k(S)^*\subset G$ (the multiplicative semigroup of $k(S)$ is
embeddable into a group), but $k(S)\nsubseteq  D$  ($k(S)$ is not
embeddable into any division ring) (see \cite{b69}).

\textbf{(b) }\ \ Group algebras $k(G)$ for Tits systems $(G,B,N,S)$
(see \cite{bou}). Here G has strong $\Gamma$-Gr\"{o}bner-Shirshov
basis, where $\Gamma=B$ and $\Gamma$-normal form is the Bruhat
normal form.

\subsection{ Tensor product of free algebras}

In A. A. Mikhalev and A. A. Zolotykh \cite{Mik3}, a CD-lemma for the
algebra $k[ X] \otimes k\langle Y\rangle$ was found, where $k[X]$ is
a polynomial algebra generated by $X$ and $k\langle X\rangle$ is a
free algebra.

In this section, we introduce the CD-lemma for tensor product
$k\langle X\rangle \otimes k\langle Y\rangle$ of free algebras,
which is from L. A. Bokut, Y. Q. Chen and Y. S. Chen \cite{cc}.

Let $X$ and $Y$ be linearly ordered sets, $S=\{yx=xy|x\in X, \ y\in
Y\}$. Then, with the deg-lex order ($y>x$ for any $x\in X, \ y\in
Y$) on $(X\cup Y)^*$, $S$ is a Gr\"{o}bner-Shirshov basis in
$k\langle X\cup Y\rangle$.  Then,  the set
$$
N=X^*Y^*=Irr(S)=\{u=u^Xu^Y| u^X\in X^* \ and \ u^Y\in Y^*\}
$$
is the normal words of the tensor product of the free algebras
$$
k\langle X\rangle \otimes k\langle Y\rangle=k\langle X\cup Y \ | \ S
\rangle.
$$

Let $kN$ be a $k$-space spanned by $N$. For any
$u=u^Xu^Y,v=v^Xv^Y\in N$, we define the multiplication of the normal
words as follows
$$
uv=u^Xv^Xu^Yv^Y\in N.
$$
Then,  $kN$ clearly coincides with the tensor product $k\langle
X\rangle \otimes k\langle Y\rangle$.

 Now, we order the set $N$. For any
$u=u^Xu^Y,v=v^Xv^Y\in N$,
$$
u>v\Leftrightarrow |u|>|v|  \ or \ (|u|=|v| \ and \ (u^X>v^X \  or \
(u^X=v^X \ and \ u^Y>v^Y))),
$$
where $|u|=|u^X|+|u^Y|$ is the length of $u$. It is obvious that $>$
is a monomial order on $N$. Such an order is also called the deg-lex
order on $N=X^*Y^*$.

Let $f$ and $g$ be monic polynomials of $kN$ and $w=w^Xw^Y\in N$.
Then we have found 16 types of compositions of inclusion and
intersection.

 $S$ is called a \emph{Gr\"{o}bner-Shirshov basis} in $kN=k\langle
X\rangle \otimes k\langle Y\rangle$ if all compositions of elements
in $S$ are trivial modulo $S$.

\begin{theorem}
 Let
$S\subseteq k\langle X\rangle \otimes k\langle Y\rangle$ with each
$s\in S$ monic and $``<"$ the deg-lex order on $N=X^*Y^*$ as before.
Then the following statements are equivalent:
\begin{enumerate}
\item[(1)]\ $S$ is a Gr\"obner-Shirshov basis in $k\langle
X\rangle \otimes k\langle Y\rangle$.
\item[(2)]\ $
f\in Id(S)\Rightarrow\overline{f}=a\overline{s}b$ for some $a,b\in N
, \ s\in S$.
\item[(3)]\
$Irr(S)=\{w \in N|w\neq a\overline{s}b, \ a,b\in N, \ s\in S\}$ is a
$k$-linear basis for the factor $k\langle X\cup Y|yx=xy, x\in X,
y\in Y\rangle/Id(S)$.
\end{enumerate}
\end{theorem}

\section{ Application of known CD-Lemmas}

\subsection{Schreier extensions of groups}

Consider a Schreier extension of group
$$
1\rightarrow A\rightarrow G \rightarrow  B\rightarrow 1.
$$

Then we have Schreier's theorem (see \cite{s1}): A group $G$ is a
Schreier extension of $A$ by $B$ if and only if there exist a factor
set $\{(b,b')|b,b'\in B\}$ of  $B$ in $A$ and $\{b: \ A\rightarrow
A, \ a\mapsto a^b \mbox{ is an automorphism}\}$ such that for any
$b,b',b''\in B, a\in A$,
\begin{eqnarray*}
(b,b')a^{bb'}=a^{[bb']}(b,b') \ \ \mbox{ and } \ \ \
(b,b'b'')(b',b'')=(bb',b'')(b,b')^{b''},
\end{eqnarray*}
where $[bb']$ is the product of elements $b,b'$ in $G$.

 M. Hall in his book \cite{h} wrote down the following statement: ``It
is difficult to determine the identities [in A] leading to
conditions for  an extension'', where the group $B$ is presented by
generators and relations.

In a  recent paper, Y. Q. Chen  \cite{c1}, by using
Gr\"{o}bner-Shirshov bases,  the structure of Schreier extensions of
groups is completely characterized and an algorithm is given to find
conditions for any Schreier extension of a group $A$ by $B$, where
$B$ is presented by a presentation. Therefore, the above problem of
M. Hall is solved.

Let $A,B$ be groups. By a factor set of $B$ in $A$, we mean a subset
of $A$ which is related to the presentation of $B$, see below.

Let the group $B$ be presented as semigroup by generators and
relations: $B=sgp\langle Y|R \rangle $, where $R$ is a
Gr\"{o}bner-Shirshov basis for $B$ with the deg-lex order $<_B$ on
$Y^*$. For the sake of  convenience, we can assume that $R$ is a
minimal Gr\"{o}bner-Shirshov basis in a sense that the leading
monomials are not contained each other as subwords, in particular,
they are pairwise different. Let $G$ be as in the following Theorem
\ref{th4.1}, where $A_1=A\backslash\{1\}, \ S=\{aa'=[aa'],
v=h_v\cdot(v), \ ay=ya^y |v\in\Omega, \  a,a'\in A_1, \ y\in Y\}$,
$\{(v)|v\in\Omega\}\subseteq A$ a factor set of $B$ in $A$, $\psi_y:
\ A\rightarrow A, \ a\mapsto a^y$ an automorphism.

We define a tower order on $(A_1\cup Y)^*$ which extends the order
$<_B$ on $Y^*$.

For $w_1=w=v_1c=dv_2, \ v_1,v_2\in \Omega, \ c,d\in Y^*, \ deg(v_1)
+ deg(v_2) > deg (w)$, we have,
$$
f_{v_1}c - df_{v_2} =dh_{v_2} - h_{v_1}c\equiv 0  \ \ mod(R,w)
$$
It means that there exists a $z\in Y^*$ such that
\begin{eqnarray*}
h_{v_1}c \equiv dh_{v_2} \equiv z \ \ mod(R,w)
\end{eqnarray*}
and thus, there exist $\xi_{(v_1,v_2)_{_w}}(v), \
\zeta_{(v_1,v_2)_{_w}}(v)\in A$ such that
\begin{eqnarray}\label{e4.1}
 g\equiv
z(\xi_{(v_1,v_2)_{_w}}(v)-\zeta_{(v_1,v_2)_{_w}}(v)) \ mod(S,w_{_1})
\end{eqnarray}
 where $\xi_{(v_1,v_2)_{_w}}(v)$ and
$\zeta_{(v_1,v_2)_{_w}}(v)$ are functions of $\{(v)|v\in\Omega\}$,
and
$g=(v_1-h_{v_{_1}}\cdot(v_1),v_2-h_{v_{_2}}\cdot(v_2))_{w_{_1}}$.

In fact, by the previous formulas, we have an algorithm to find the
functions $\xi_{(v_1,v_2)_{_w}}(v)$ and $
\zeta_{(v_1,v_2)_{_w}}(v)$.

\begin{theorem}(\cite{c1})\label{th4.1}
Let  $A,B$ be groups, $B=sgp\langle Y|R \rangle$, where
$R=\{v-h_v|v\in\Omega\}$ is a minimal Gr\"{o}bner-Shirshov basis for
$B$ and $v$ the leading term of the polynomial $f_v=v-h_v\in R$. Let
$$
G=E(A,Y,a^y,(v))=sgp\langle A_1\cup Y| S\rangle
$$
where $A_1=A\backslash\{1\}, \ S=\{aa'=[aa'], \ v=h_v\cdot(v), \
ay=ya^y |v\in\Omega, \ a,a'\in A_1, \ y\in Y\}$, $\psi_y: \
A\rightarrow A, \ a\mapsto a^y$ an automorphism,
$\{(v)|v\in\Omega\}\subseteq A$ a factor set of $B$ in $A$ with
$(v)=1$ if $f_v=y^{\epsilon}y^{-\epsilon}-1, \ y\in Y, \
\epsilon=\pm 1$.
\begin{enumerate}
\item[(i)]
For the tower order, $S$ is a Gr\"{o}bner-Shirshov basis for $G$ if
and only if for any $v\in\Omega, \ a\in A$ and any composition
$(f_{v_1},f_{v_2})_{_w}$ of $R$ in $k\langle Y \rangle$,
\begin{eqnarray}\label{e4.2}
(v) a^{v}= a^{h_v}(v) \ \ \mbox{ and } \ \ \
\xi_{(v_1,v_2)_{_w}}(v)=\zeta_{(v_1,v_2)_{_w}}(v)
\end{eqnarray}
hold in $A$, where $\xi_{(v_1,v_2)_{_w}}(v)$ and
$\zeta_{(v_1,v_2)_{_w}}(v)$ are defined by (\ref{e4.1}). Moreover,
if this is the case, $G$  is a Schreier extension of $A$ by $B$ in a
natural way.
\item[(ii)]
A group $C$ is a Schreier extension of $A$ by $B$ if and only if
there exist $\{a^y|y\in Y, \ A\rightarrow A, \ a\mapsto a^y \mbox{
is an isomorphism}\}$ and a factor set $\{(v)|v\in\Omega\}$ of  $B$
in $A$  with $(v)=1$ if $f_v=y^{\epsilon}y^{-\epsilon}-1, \ y\in Y,
\ \epsilon=\pm 1$
 such that (\ref{e4.2}) holds. Moreover, if this is the
case, $C\cong G=E(A,Y,a^y,(v))=sgp\langle A_1\cup Y| S\rangle$.
\end{enumerate}
\end{theorem}

\noindent {\bf Remark.} \ In the above theorem, let the group
$A=gp\langle X|R_A \rangle$ be also presented by generators and
relations, where $R_A=\{u=f_u|u\in\Lambda\}$ is a
Gr\"{o}bner-Shirshov basis for $A$ with the deg-lex order $<_A$ on
$X^*$. Then, by replacing $A_1$ with $X$, $aa'=[aa']$ with $R_A$ and
$x$ with $a$, the results hold.

 As a corollary of the above theorem, in \cite{c1}, by using the result
 in Y. Q. Chen and C. Y. Zhong in \cite{cz1},  a  criteria in
the case that $B$ is HNN-extension is formulated.

Another solution to the M. Hall's problem can be referred in
\cite{pride}.

\subsection{Extensions of algebras }

In the paper of Y. Q. Chen \cite{c2}, he gave the same kind of
answer to an analogy of the M. Hall's problem in the above section
for Schreier extensions of algebras.

\begin{definition} Let $k$ be a field, $M,B, \cal{R}$ $k$-algebras (not necessarily with 1).
Then $\cal{R}$ is called an extension of $M$ by $B$ if  $M^2=0$,
where $M$ is an ideal of $\cal{R}$ and $\cal{R}$$/M\cong B$ as
algebras. Such an extension is called a singular extension in
\cite{ho}.
\end{definition}

The following classical result is known. Let $M,B, \cal{R}$ be
$k$-algebras with $M^2=0$. Then $\cal{R}$ is an extension of  $M$ by
$B$ if and only if $M$ is a $B$-bimodule and there exists a factor
set $\{(b,b')|b,b'\in B\}$ of $B$ in $M$ such that for any
$b,b',b''\in B$,
$$
b(b',b'')-(bb',b'')+(b,b'b'')-(b,b')b''=0.
$$

In \cite{c2}, by using Gr\"{o}bner-Shirshov bases,  the structure of
extensions of algebras is completely characterized and an algorithm
is given to find conditions for any  extension of an algebra $M$ by
$B$, where $B$ is presented by a presentation.

As results, by using this theorem, in \cite{c2}, a characterization
theorem of the extension of $M$ by $B$ is given, when  $B$ is a
cyclic algebra, free commutative algebra, universal envelope of a
Lie algebra, and Grassmann algebra,
 respectively.

\subsection{Anti-commutative algebras}

In the paper of L. A. Bokut, Y. Q. Chen and Y. Li \cite{bcly}, an
application of Shirshov's CD-lemma for anti-commutative algebras
(see \cite{S3}) is given. This application gives an anti-commutative
Gr\"{o}bner-Shirshov basis of a free Lie algebra.

Let $X=\{x_i|i\in I\}$ be a linear ordered set. Let $X^{\ast \ast }$
be the set of all non-associative words $(u)$ in $X$. We assume that
$(u)$ is a bracketing of $u$. Then we define normal words
$N=\{[u]\}$ and order them by using induction on the length
$n=|[u]|$ of $[u]$:
\begin{description}
\item If $n=1$, then $[u]=x_{i}$ is a normal word. Define $%
x_{i}>x_{j}$ if $i>j$.

Let $N_{n-1}=\{[u]|[u]\mbox{ is  a  normal  word  and }|[u]|\leq
{n-1}\}$, $n>1$ and suppose that $``<"$ is a well order on
$N_{n-1}$. Then
\item  If $n>1$ and $(u)=((v)(w))$ is a word of length $n$, then $(u)$
is a normal word, if and only if
\begin{description}
\item  both $(v)$ and $(w)$ are normal words, that is, $(v)=[v]$ and $(w)=[w]$, and
\item  $[v]>[w]$.
\end{description}
\end{description}
Let $[u]$, $[v]$ be normal words of length $\leq n$. Then $[u]<[v]$,
if and only if one of the following three cases holds:
\begin{description}
\item $|[u]|<n$, $|[v]|<n$ and $[u]<[v]$.
\item  $|[u]|<n$ and $|[v]|=n$.
\item  If $|[u]|=|[v]|=n$, $[u]=[[u_{1}][u_{2}]]$ and $[v]=[[v_{1}][v_{2}]]$, then $%
[u_{1}]<[v_{1}]$ or $([u_{1}]=[v_{1}]\ and\ [u_{2}]<[v_{2}])$.
\end{description}
It is clear that the order ``$<$" on $N$ is a well order.

Let $AC(X)$ be a $k$-space spanned by normal words. Now, we define
the product of normal words by the following way:
\begin{equation*}
\widetilde{[u][v]}=\left\{
\begin{array}{r@{\quad:\quad}l}
[[u][v]] & [u]>[v] \\
-[[v][u]] & [u]<[v] \\
0\text{ \ \ } & [u]=[v]%
\end{array}%
\right.
\end{equation*}

Then $AC(X)$ is the free anti-commutative algebra generated by $X$.

Let $S\subset AC(X)$ be a set of monic polynomials, $s\in S$ and
$(u)\in X^{**}$. We define $S$-word $(u)_s$ by induction as a
non-associative word in $X\cup S$ with only one occurrence of $s\in
S$. An $S$-word $(u)_s$ is called a normal $S$-word if $(u)_{\bar s
}=(a[\bar s ]b)$ is a normal word.

There is only one kind of compositions that is inclusion one.

\begin{theorem}( \cite{S3}, \cite{bcly})
Let $S\subset AC(X)$ be a nonempty set of monic polynomials and the
order $``<"$ as before. Then the following statements are
equivalent:
\begin{enumerate}
\item [(i)] $S$ is a Gr\"{o}bner-Shirshov basis in $AC(X)$.

\item [(ii)] $f\in Id(S)\Rightarrow [\bar f] =[a[\bar s ]b]$ for some $s\in S\
and\ a,b\in X^*$, where $[as b]$ is  normal $S$-word.

\item [(iii)] $Irr(S)=\{[u]\in N |[u]\ne [a[\bar s] b]\ a,b\in X^*,\ s\in S \mbox{ and }
[as b] \mbox{ is a normal } S\mbox{-word}\}$ is a basis of the
algebra $AC(X|S)$.
\end{enumerate}
\end{theorem}

By using this theorem, a Gr\"{o}bner-Shirshov basis $S$ in $AC(X)$
is given in \cite{bcly}  which shows that the Hall words in $X$
forms a basis for the free Lie algebra $Lie(X)$, where
$S=\{([u][v])[w]-([u][w])[v]]-[u]([v][w]) \ |[u]>[v]>[w]\ and \
[u],[v],[w]\ are\ Hall\ words\}$.

\subsection{Akivis algebras }

This section is from the paper of Y. Q. Chen and Y. Li \cite{cl}.

An Akivis algebra is a vector space $A$ over a field $k$ endowed
with a skew-symmetric bilinear product $[x,y]$ and a trilinear
product (x,y,z) that satisfy the identity
$[[x,y],z]+[[y,z],x]+[[z,x],y]=(x,y,z)+(z,x,y)+(y,z,x)-(x,z,y)-(y,x,z)-(z,y,x)$.
For any (non-associative) algebra $B$, one may obtain an Akivis
algebra $Ak(B)$ by considering in $B$ the usual commutator
$[x,y]=xy-yx$ and associator $(x,y,z)=(xy)z-x(yz)$.

The CD-lemma for non-associative algebras is invented by Shirshov in
\cite{S3}. By applying this lemma, in \cite{cl}, a
Gr\"{o}bner-Shirshov basis in $AC(X)$ is given for the universal
enveloping algebra of an Akivis algebra which gives an another proof
of I. P. Shestakov's result (see \cite{She}) that any Akivis algebra
is linear. An Akivis algebra $A$ is linear if $A$ can be embedded in
some non-associative algebra $B$ with above operations.

\begin{theorem} (\cite{cl})
Let $(A,+,[,],(,,))$ be an Akivis algebra with a linearly ordered
basis $\{e_i|\ i\in I\}$. Let
$$
[e_i,e_j]=\sum\limits_{m}\alpha_{ij}^me_m,\
(e_i,e_j,e_k)=\sum\limits_{n}\beta_{ijk}^ne_n,
$$
where $\alpha_{ij}^m,\beta_{ijk}^n\in k$. We denote
$\sum\limits_{m}\alpha_{ij}^me_m$ and
$\sum\limits_{n}\beta_{ijk}^ne_n$ by $\{e_ie_j\}$ and
$\{e_ie_je_k\}$, respectively. Let
$$
U(A)=M( \{e_i\}_I| \ e_ie_j-e_je_i=\{e_ie_j\}, \
(e_ie_j)e_k-e_i(e_je_k)=\{e_ie_je_k\}, \ i,j,k \in I)
$$
be the universal enveloping algebra of $A$. Let
\begin{eqnarray*}
S&=&\{e_ie_j-e_je_i-\{e_ie_j\} \ (i>j), \
(e_ie_j)e_k-e_i(e_je_k)-\{e_ie_je_k\} \ (i,j,k\in I),
\\
&&e_i(e_je_k)-e_j(e_ie_k)-\{e_ie_j\} e_k-\{e_je_ie_k\}+\{e_ie_je_k\}
\ (i>j, \ k\geq j) \}.
\end{eqnarray*}
Then \begin{enumerate}
\item[(i)] \ $S$ is a Gr\"{o}bner-Shirshov basis for $U(A)$,
\item[(ii)] \ $A$ can be embedded into the universal
enveloping algebra $U(A)$.
\end{enumerate}
\end{theorem}

\subsection{Some one-relator groups }

This section contains the results of Y. Q. Chen and C. Y. Zhong in
\cite{cz2}.

It is well known the Magnus algorithm for a solution of the word
problem for any one-relator group.

Using the Magnus rewriting procedure (see R. C. Lyndon and P. E.
Schupp \cite{ly}), one may embed any one-relator group into a tower
of HNN-extensions. For towers of HNN-extensions of groups, L. Bokut
(see \cite{bku}) developed a method of groups with the standard
normal forms in 1965. Actually, for these groups,
 Gr\"obner-Shirshov bases are also ``standard" in a sense, so we may speak about
``groups with the standard  Gr\"obner-Shirshov bases".

By the way, we give here a story about applications of groups with
the standard normal forms. At our seminar, we have studied the
following interesting paper

 \textit{K. Kalorkoti, Turing degree and the word and
conjugacy problem for finitely presented groups, Internet.}

 Actually, it is a part of his Thesis in London
University, 1979.

C. Y. Zhong was the speaker for almost two months at the seminar and
she pointed out that K. Kalorkoti used successfully the method of L.
A. Bokut on  groups with the standard normal forms.

As a result of this study, we suggested K. Kalorkoti to publish his
paper in the SEA Bull. Math.  (see \cite{ka}).

 Now, we back to one-relator groups. There is a chance
that any tower of HNN-extensions produced by the Magnus method has
the standard Gr\"obner-Shirshov basis and the standard normal form.
In particular, it would give another algorithm for the word problem
for any one-relator group.

 The problem is that the Magnus embedding is not easy
to observe and write down explicitly. Hence, we need to go step by
step. Any one-relator group can be effectively embedded into
one-relator group with two generators
$$
gp\langle X|r \rangle\hookrightarrow G=gp\langle x,\
y|x^{n_1}y^{m_1}\cdots x^{n_k}y^{m_k}=1 \rangle,
$$
where $x_i\mapsto x^{-q}yx^l$ (W. Magnus, see also \cite{ly}), $n_i,
\ m_i\neq0,  \ k\geq0$.
 We call $k$ the depth of $G$.
\begin{theorem}(\cite{cz2})
Any two-generator one-relator group of the depth 3 is effectively
Magnus embeddable into a tower of HNN-extensions, which is a group
with the effective standard Gr\"obner-Shirshov basis and effective
standard normal form.
\end{theorem}

\subsection{The Chinese monoid}

The Chinese monoid $CH(X)$ on a well ordered set $X$ has the
following  defining relations:
$$
cba=bca=cab, \ \ c\geq b\geq a, \ \ a,b,c\in X.
$$
A fundamental paper on the Chinese monoid has been published in 2001
\cite{cek}.

In the paper of Y. Q. Chen and J. J. Qiu \cite{cq2}, a
Gr\"obner-Shirshov basis for the Chinese monoid is found.

\begin{theorem} (\cite{cq2}) Let $S=\{x_ix_jx_k-x_jx_ix_k, \
x_ix_kx_j-x_jx_ix_k, \  x_ix_jx_j-x_jx_ix_j, \ x_ix_ix_j-x_ix_jx_i,
\ x_ix_jx_ix_k-x_ix_kx_ix_j, \ x_i, x_j, x_k \in X,\ i>j>k\}$. Then
\begin{enumerate}
\item[(i)] $sgp\langle X|T\rangle=sgp\langle X|S\rangle$ and with deg-lex order,
$S$ is a Gr\"{o}bner-Shirshov basis of the Chinese monoid $CH(X)$.
\item[(ii)]The $Irr(S)$ is the set which consists of words on $X$ of the
form $u_n=w_1w_2\cdots w_n, n\geq0$, where
\begin{eqnarray*}
w_1 &=& x_1^{t_{11}}\\
w_2 &=& (x_2x_1)^{t_{21}}x_2^{t_{22}}\\
w_3 &=&  (x_3x_1)^{t_{31}}(x_3x_2)^{t_{32}}x_3^{t_{33}}\\
& & \hspace{0.4cm} \cdots\\
w_n &=&
(x_nx_1)^{t_{n1}}(x_nx_2)^{t_{n2}}\cdots(x_nx_{n-1})^{t_{n(n-1)}}
x_n^{t_{nn}}
\end{eqnarray*}
with $x_i\in X $, $x_1<x_2< \cdots <x_n$ and all exponents are
non-negative integers.
\end{enumerate}
\end{theorem}

 Then $Irr(S)$ coincides  with the set of staircase
words of the  paper by J. Cassaigne et al. \cite{cek}. Also, the
insertion algorithm of J. Cassaigne et al. \cite{cek} coincides with
the elimination of leading words algorithm.

\subsection{ Markov and Artin normal form theorem for braid groups}

In the paper of L. A. Bokut, V. V. Chaynikov and K. P. Shum
\cite{bcs}, the authors present the classical results of
Artin-Markov on braid groups by using the Gr\"obner-Shirshov bases.
As an application, one can reobtain the normal form of
Artin-Makov-Ivanovskiy as an easy corollary.

\section{ Expository papers}

\subsection{ Gr\"obner and Gr\"obner-Shirshov bases: an elementary
approach }

There is an elementary approach in L. A. Bokut and K.P. Shum
\cite{bs2} to Gr\"{o}bner-Shirshov bases theory with quite  a few
examples, including the example for Lie algebras.

\subsection{Shirshov's CD-lemma for Lie
algebras}

What is now called the Gr\"{o}bner-Shirshov method for Lie
 algebras originally invented by A. I. Shirshov in 1962 (\cite{Sh}). Actually, that
 paper based on his paper \cite{S58} where Shirshov invented a new
 linear basis for a free Lie algebra now called the  Lyndon-Shirshov
 basis (it was defined independently in the paper \cite{C} in the same year).
Remarkably,  Lyndon--Shirshov basis is a particular case of the
 series of bases of a free Lie algebra invented by A. I. Shirshov in his
 Candidate of Doctor of Science Thesis (Moscow State University, 1953, advisor was A. G. Kurosh)
 and published in 1962 (\cite{S4})
(cf. \cite{R} where these bases are called Hall Bases).
 We now cite the Zbl review by P. M. Cohn \cite{Cohn}] of the paper \cite{S58}:
  ``The author varies the usual construction of basis commutators
 in Lie rings by ordering words lexicographically and not by length.
 This is used to give a very short proof of the theorem
 (Magnus \cite{M}, Witt \cite{W})
 that the Lie algebra obtained from a free associative algebra is
 free. Secondly he derives Friedrich's criterion (this Zbl 52,45)
 for Lie elements. As the third application he proves that every Lie
 algebra $L$ can be embedded in a Lie algebra $M$ such that in $M$
 any subalgebra of countable dimension is contained in a 2-generated
 subalgebra." We would like to add that it was a beginning
 of Gr\"{o}bner-Shirshov bases theory for Lie and associative
 algebras. Lemma 4 of the paper, on special bracketing of a regular
 (Lyndon-Shirshov) associative word with a fix regular subword,
 leads to the algorithm of elimination of the leading
word of one Lie polynomial in other Lie polynomial , i.e., to the
reduction procedure that is very familiar in the cases of
associative and associative-commutative polynomials. Also this Lemma
4 leads to the crucial notion of composition of two Lie polynomials
that will be defined lately in the paper \cite{Sh}. Last but not
least, Shirshov \cite{S53} proved the following result for
connections of some ideals of free Lie and free associative
algebras.

Let $Lie(a,b)$ be the Lie algebra of Lie polynomials of $k\langle
a,b\rangle$ (it is the free Lie algebra over a set $\{a,b\}$ and a
field $k$). Let $J=J([a^2b^kab]=[[a[[ab]\cdots b]]][ab],\ k\geq 1)$
be the Lie ideal of $Lie(a,b)$ generated by $\{[a^2b^kab],\ k\geq
1\}$ and $I$ the associative ideal of $k\langle a,b\rangle$
generated by $J$. Then, $I\cap Lie(a,b)=J$. The proof is dealing
with leading monomials of Lie and associative polynomials. This
result and its proof are the real beginning of Gr\"{o}bner-Shirshov
bases theory for Lie and associative algebras.

As for the  paper \cite{Sh} itself, it is a fully pioneer paper in
the
 subject. He defines a notion of the composition $(f,g)_w$ of two
 Lie (associative) polynomials relative to an associative word $w$
 (it was called lately by $S$-polynomial for commutative polynomials by
 B. Buchberger \cite{bu65} and \cite{bu70}).
 It leads to the algorithm for construction of a Gr\"{o}bner-Shirshov
 basis ($GSB(S)$) of Lie (associative) ideal generated by some set
 $S$: to joint to S all nontrivial compositions and to eliminate
leading monomials of one polynomial of S in  others. Shirshov proves
the lemma, now known as the Composition, or
 Composition-Diamond lemma, that if $f \in Id_{Lie}(S)$, then
 $\overline{f}$, the leading associative word of $f$, has a form $\overline{f}=u
 \overline{s}v$, where $s \in GSB(S), \ u,v \in X^*$. Several years
 later the first author formulated this lemma in the modern form (see
\cite{b72}). Let $S$ be a set of Lie polynomials that is complete
under
 composition ( i.e., any composition of polynomials of $S$ is trivial;
  on the other word, $S$ is a Gr\"{o}bner-Shirshov
 basis). Then if $f \in Id_{Lie}(S)$, then $\overline{f}=u
 \overline{s}v$, where $s \in S, \ u,v \in X^*$. Of course, from
 Shirshov's Composition-Diamond lemma it follows that the set
 $Irr(S)$ of $S-$irreducible Lyndon-Shirshov words constitutes a linear
 basis of the quotient algebra $Lie(X)/ Id(S)$. The converse is
 also true.

Lately explicitly Shirshov's Composition-Diamond lemma for
associative algebra was formulated by L. A. Bokut \cite{b76} in 1976
and G. Bergman \cite{b} in 1978.

The paper by L. A. Bokut and Y. Q. Chen \cite{bc} contains a
comprehensive proof of the Shirshov's Composition-Diamond lemma for
Lie algebras. We follow Shirshov's ideas of his papers \cite{S53},
\cite{S58}, \cite{S4} and \cite{Sh}. In particular, we prove all
necessary properties of both associative (see also \cite{lyn}) and
non-associative Lyndon-Shirshov words by using the Shirshov's
elimination in \cite{S53} (that is, the so called  Lazard
elimination in \cite{L} and \cite{R}).

\subsection{CD-lemma for modules}

Composition-Diamond lemma for modules was first formulated and
proved by S.-J. Kang and K.-H. Lee in \cite{kl1} and \cite{kl2}.
According to their approach, a Gr\"{o}bner-Shirshov basis of a
cyclic module $M$ over an algebra $A$ is a pair $(S,T)$, where $S$
is the set of the defining relations of $A, \ A=k\langle
X|S\rangle$, and $T$ is the defining relations for the $A$-module
$_A M= _A\!M(e|T)$. Then Kang-Lee's Lemma says that $(S,T)$ is a
Gr\"{o}bner-Shirshov pair for the $A$-module $_A M=_A M(e|T)$ if $S$
is a Gr\"{o}bner-Shirshov basis of $A$ and $T$ is closed under the
right-justified composition with respect to $S$, and for $f\in S, \
g\in T$, such that $(f,g)_w$ is defined and $(f,g)_w\equiv 0 \
mod(S,T,w)$. A  composition  $(f,g)_w$ is called right-justified if
$w=\overline{f}=a\overline{g}$ for some $a\in X^*$.

Some years later, E. S. Chibrikov \cite{ch04} suggested a new
Composition-Diamond lemma for modules that treat any module as a
factor module of ``double-free" module.

Let $X,Y$ be sets and $mod_{k\langle X\rangle} \langle Y\rangle$ a
free left $k\langle X\rangle$-module with the basis $Y$. Then
$mod_{k\langle X\rangle} \langle Y\rangle=\oplus_{y\in Y} k\langle
X\rangle y$ is called a ``double-free"  module.

Suppose that $<$ is a monomial order on $X^*$, $<$ a well order on
$Y$  and $X^*Y=\{uy|u\in X^*, \ y\in Y\}$. We define an order
$``\prec"$ on $X^*Y$: for any $w_1=u_1y_1,w_2=u_2y_2\in X^*Y$,
\begin{equation}\label{e5.3}
w_1\prec w_2\Leftrightarrow u_1<u_2 \ \ \mbox{ or } u_1=u_2, \
y_1<y_2
\end{equation}

 Let $S\subset
mod_{k\langle X\rangle} \langle Y\rangle$ with each $s\in S$ monic.
We define a composition $(f,g)_w=f-ag$, where $w=\bar{f}=a\bar{g}, \
a\in X^*, \ f,g\in mod_{k\langle X\rangle} \langle Y\rangle$ are
monic.

If $(f,g)_w=f-ag=\sum \alpha_i a_i s_i$, where $\alpha_i\in k, \ a_i
\in X^*, \ s_i \in S$ and $a_i \overline{s}_i\prec w$, then this
composition is called trivial modulo $(S,w)$ and is denoted  by $
(f,g)_w\equiv 0 \ mod(S,w). $

\begin{definition}\label{d3}(\cite{ch04})
Let $S\subset mod_{k\langle X\rangle} \langle Y\rangle$ be a
non-empty set  with each $s\in S$ monic. Let the order ``$\prec$" be
as before. Then we call $S$ a Gr\"obner-Shirshov basis in the module
$mod_{k\langle X\rangle} \langle Y\rangle$ if all the compositions
of polynomials in $S$ are trivial modulo $S$.
\end{definition}

\begin{lemma}\label{l5.3.2}
(\cite{ch04}, Composition-Diamond lemma for ``double-free"
modules)\label{l2} Let $S\subset mod_{k\langle X\rangle} \langle
Y\rangle$ be a non-empty set  with each $s\in S$ monic and
$``\prec"$ the order on $X^*Y$ as before. Then the following
statements are equivalent:
\begin{enumerate}
\item[(i)]\ $S$ is a Gr\"obner-Shirshov basis in $mod_{k\langle X\rangle} \langle Y\rangle$.
 \item[(ii)]\  $
f\in k\langle X\rangle S\Rightarrow\overline{f}=a\overline{s}$ for
some $a\in X^*, \ s\in S$.
\item[(iii)]\
$Irr(S)=\{w \in X^*Y|w\neq a\overline{s}, \ a\in X^*, \ s\in S\}$ is
a $k$-linear basis for the factor $mod_{k\langle X\rangle} \langle
Y|S\rangle=mod_{k\langle X\rangle} \langle Y\rangle /k\langle
X\rangle S$.
\end{enumerate}
\end{lemma}

As applications of Lemma \ref{l5.3.2},  Y. Q. Chen, Y. S. Chen and
C. Y. Zhong \cite{ccz},  found the  Gr\"obner-Shirshov bases for
highest weight modules over  Lie algebra $sl_2$, Verma modules over
Kac-Moody algebras, Verma modules over Lie algebras of coefficients
of free conformal algebras and the universal enveloping modules for
 Sabinin algebras. (The last modules are defined in \cite{pe}).

\end{document}